\documentclass[12pt]{article}

\makeatletter
\newfont{\footsc}{cmcsc10 at 8truept}
\newfont{\footbf}{cmbx10 at 8truept}
\newfont{\footrm}{cmr10 at 10truept}
\renewcommand{\ps@plain}{%
\renewcommand{\@oddfoot}{\footsc the electronic journal of combinatorics
  {\footbf 9} (2002), \#R39\hfil\footrm\thepage}}
\makeatother
\usepackage{amsmath}
\usepackage{amssymb}
\usepackage{amsthm}
\newtheorem*{thm}{Theorem}
\newtheorem*{conj}{Conjecture}
\theoremstyle{definition}
\newtheorem*{defn}{Definition}
\newtheorem*{lemma}{Lemma}
\newcommand{\by}{\times}
\newcommand{\lambdac}{{\lambda^{\mathrm c}}}
\newcommand{\muc}{{\mu^{\mathrm c}}}
\newcommand{\pa}[1]{{\langle{#1}\rangle}}
\newcommand{\ww}{{\mathbf w}}
\newcommand{\hh}{{\mathtt h}}
\newcommand{\vv}{{\mathtt v}}
\newcommand{\hl}[3]{\put(#1,#2){\line(1,0){#3}}}
\newcommand{\vl}[3]{\put(#1,#2){\line(0,1){#3}}}
\newcommand{\yd}[2]{\,\,{\begin{picture}#1{#2}\end{picture}}\,\,}
\newcommand{\mybox}{{\begin{picture}(1,1) \hl001 \hl0{-1}1 \vl0{-1}1 \vl1{-1}1 \end{picture}}}
\newcommand{\row}[2]{\multiput(0,#1)(1,0){#2}{\mybox}}
\newcommand{\rb}[1]{\raisebox{#1\unitlength}}
\newcommand{\rr}[2]{\rb{-1}{\yd{(#1,2)}{\row2{#1}\row1{#2}}}}  
\newcommand{\rrr}[3]{\rb{-2}{{\yd{(#1,3)}{\row3{#1}\row2{#2}\row1{#3}}}}} 

\newcommand{\rrrrr}[5]{\rb{-4}{{\yd{(#1,5)}{\row5{#1}\row4{#2}\row3{#3}\row2{#4}\row1{#5}}}}}
\newcommand{\pad}{\hspace{\unitlength}}
\newcommand{\inline}[1]{{\setlength{\unitlength}{3.5pt}\rb{1.5}{#1}}}  
\newcommand{\Inline}[1]{{\setlength{\unitlength}{5.0pt}\rb{.4}{#1}}}  
\newcommand{\Srow}[2]{\multiput(-#1,#1)(1,0){#2}{\mybox}}
\newcommand{\Srr}[2]{\rb{-1}{\yd{(#1,2)(-2,0)}{\Srow2{#1}\Srow1{#2}}}}  
\newcommand{\Srrr}[3]{\rb{-2}{{\yd{(#1,3)(-3,0)}{\Srow3{#1}\Srow2{#2}\Srow1{#3}}}}} 

\title{Linearly Independent Products of \\ Rectangularly Complementary \\ Schur Functions}
\author{Michael Kleber \\
\small Department of Mathematics \\[-0.8ex]
\small Brandeis University\\[-0.8ex]
\small Waltham, MA 02454\\[-0.8ex]
\small \texttt{kleber@brandeis.edu}}

\date{\small Submitted: 11 September 2002; Accepted: 16 September 2002.\\
\small MR Subject Classification: 05E05}

\begin{document}

\maketitle

\begin{abstract} \setlength{\parindent}{0pt}
Fix a rectangular Young diagram $R$, and consider all the products of Schur functions $s_\lambda s_\lambdac$, where $\lambda$ and $\lambdac$ run over all (unordered) pairs of partitions which are complementary with respect to $R$.

{\bf Theorem:} The self-complementary products, $s_\lambda^2$ where $\lambda=\lambdac$, are linearly independent of all other $s_\lambda s_\lambdac$.

{\bf Conjecture:} The products $s_\lambda s_\lambdac$ are all linearly independent.
\end{abstract}

\section{Introduction}

The Schur functions $s_\lambda$, where $\lambda$ is a partition or its graphical representation as a Young diagram, play a pivotal role in the theory of symmetric functions and its many applications.  The literature abounds with linear identities among products of two Schur functions, many of which follow from expressing $s_\lambda$ as a determinant of a certain matrix, to which the Pl\"ucker relations are readily applied.  This paper instead presents a linear {\em independence} result among some products of pairs of Schur functions.

Schur functions are indexed by partitions $\lambda = \pa{\lambda_1,\ldots,\lambda_r}$ of all integers $n\geq0$; we order the parts to be weakly decreasing and consider the partition to end with an infinite tail of 0's, which we are free to omit.  We identify $\lambda$ with its Young diagram, stacked rows of $\lambda_1,\ldots,\lambda_r$ boxes, and write $\lambda\subseteq\mu$ if the diagram of $\lambda$ sits inside that of $\mu$, i.e. if $\lambda_i \leq \mu_i$ for all $i$.  If we fix an $a\by b$ rectangular Young diagram $R = \pa{b,\ldots,b}=\pa{b^a}$, then for any $\lambda\subseteq R$, the boxes in $R$ which remain when $\lambda$ is removed from the upper left corner can be rotated $180^\circ$ to get a new partition $\lambdac$, the complement of $\lambda$ with respect to $R$.  For example, if $R=\pa{4^3}=\inline{\rrr444}$ and $\lambda=\pa{421}=\inline{\rrr421}$ then $\lambdac=\pa{32}=\inline{\rr32}$.

We work in the ring $\Lambda$ of symmetric formal polynomials in infinitely many variables $x_1$, $x_2$,\ldots, of which the Schur functions $\{s_\lambda\}$ form a linear basis.  The structure constants of the basis, those $c_{\lambda\mu}^\nu$ such that $s_\lambda s_\mu = \sum c_{\lambda\mu}^\nu s_\nu$, are nonnegative integers called the Littlewood--Richardson numbers.  For example, in $\Lambda$:
\begin{equation}
\label{ex22}
\begin{array}{lcl}
\rr20 \rr20 &=& \rr22 + \rr31 + \rrr400 \\ 
\rr11 \pad \rr11 &=& \rr22 + \rrr211 \pad + \rrrrr11110 \\ 
\rr21 \rr10 &=& \rr22 + \rr31+\rrr211 \\
\rr22 \rr00 &=& \rr22\\ 
\end{array}
\end{equation}
These are precisely the four distinct products $s_\lambda s_\lambdac$ where $R$ is the $2\by2$ rectangle.  In the fourth line $\lambda=\pa{22}$ and $\lambdac$ is the zero partition $\pa{0}$, whose diagram is empty and whose Schur function $s_\pa{0}$ is $1\in\Lambda$.

Note that the products in~(\ref{ex22}) are easily seen to be linearly independent.  Each of the first two contains a basis element ($s_\pa{4}$ and $s_\pa{1111}$, respectively) which appears nowhere else.  We immediately know that these products cannot show up in any linear dependence; in this situation we say that some Schur function is a {\em witness} to the linear independence of the product in which it appears.  Once the first two products are dropped from consideration, the third has two witnesses to its linear independence from all remaining products; once the third is dropped, the fourth product has a witness as well.

The theorem will be proved by constructing a witness for each self-complementary $\lambda$, as in the above example (if the rectangle's dimensions are both odd, then for $\lambda$ and $\lambdac$ that differ only in which contains the center box of $R$).  The strong form of the conjecture is that the above procedure of successive elimination via witnesses always proves the linear independence of all the products.

\subsection*{Acknowledgements}
This work was inspired by Mark Adler, who in January 2001 asked me whether some statement along the lines of the theorem might be true.  His interest was in unique representability of certain special elements, the ``$\tau$-functions'',  of the vector spaces $\bigoplus_{\lambda=\lambdac} {\mathbb C} s_\lambda$.  These $\tau$-functions are solutions to the finite Pfaff lattice and relate in some special cases to Jack polynomials.  For context see~\cite{adler}, though in the end the theorem was not needed in general. 

I posed the conjecture as an open problem during the ``Symmetric Functions 2001: Surveys of Developments and Perspectives'' NATO ASI in June--July 2001, a workshop that was part of the Isaac Newton Institute's programme ``Symmetric Functions and Macdonald Polynomials;'' my thanks to the organizers for inviting me to participate.  Thanks also go to Anders Buch for his Littlewood--Richardson Calculator and to Sara Billey for a hefty chunk of computer time to run it.  The work prior to July 2001 (at MIT) was partially supported by an NSF Postdoctoral Research Fellowship.

\section{Theorem}

Fix positive integers $a$ and $b$, and let $R$ be the partition $\pa{b^a}$, whose Young diagram is an $a\by b$ rectangle.  We take all complements with respect to $R$, and say $\lambda\subseteq R$ (and $s_\lambda$) is self-complementary if $\lambda=\lambdac$.  If $a$ and $b$ are odd, there are no self-complementary $\lambda$.  In this case we will say $\lambda$ is {\em almost self-complementary} if $\lambda$ and $\lambdac$ differ only in which one contains the central box of $R$.

\begin{thm}
The squares of the self-complementary Schur functions, $s_\lambda^2$ for $\lambda=\lambdac$, are linearly independent of all other products $s_\mu s_\muc$ for $\mu\subseteq R$.
If the dimensions of $R$ are both odd, the products $s_\lambda s_\lambdac$ for $\lambda$ almost self-complementary are linearly independent.
 \end{thm}

By way of proof, we show that for each such $\lambda$ there is another partition $\pi$ which acts as a witness to the linear independence: $s_\pi$ appears in the Schur function expansion of $s_\mu s_\muc$ if and only if $\mu=\muc=\lambda$.  First we define a new family of notations for partitions.

\begin{defn}
Let $c$ and $d$ be nonnegative integers, and let $\ww$ be a word of length $c+d$ consisting of the letters $\hh$ and $\vv$ appearing $c$ and $d$ times, respectively.  Let $\lambda$ be a partition with at most $c$ parts larger than $d$.
\begin{enumerate}
\item
The {\em $\ww$-decomposition} for the Young diagram of $\lambda$ partitions the diagram into $c+d$ pieces.  The $i$th piece is part of one horizontal row or vertical column, depending on whether the $i$th letter of $\ww$ is an $\hh$ or a $\vv$, respectively.  It begins with the upper-left-most box not contained in any previous piece, and extends right or down, respectively, to the end of the diagram.
\item
The {\em $\ww$-notation} for $\lambda$ is the $c+d$-tuple of integers $\lambda^\ww = (\lambda_1^\ww,\ldots,\lambda_{c+d}^\ww)$, where $\lambda_i^\ww$ is the number of boxes in the $i$th part of the $\ww$-decomposition.
\item
The {\em $\ww$-ordering} on partitions is the lexicographic ordering on their $\ww$-notations.
\end{enumerate}
\end{defn}

For example, Figure~\ref{fig_w}(a) depicts the $\ww$-decomposition of $\mu=\pa{5411}$ with  $\ww=\hh\vv\vv\hh$;  we see that $\mu^\ww=(5,3,1,2)$.  The $\ww$-ordering has the following useful property:

\begin{lemma}
Fix a word $\ww$ as above, and take two partitions with
$\ww$-notations $\mu^\ww=(\mu_1^\ww,\ldots,\mu_{c+d}^\ww)$ and
$\nu^\ww=(\nu_1^\ww,\ldots,\nu_{c+d}^\ww)$.  Let $\pi$ be the partition
with $\pi_i^\ww = \mu_i^\ww + \nu_i^\ww$. Then $s_\pi$ appears in the
Schur expansion of $s_\mu s_\nu$, and it is the maximal partition (in the
$\ww$-ordering) to appear.
\end{lemma}

\proof
We show that $s_\pi$ appears in $s_\mu s_\nu$ by presenting a Littlewood--Richardson tableau of shape $\pi/\mu$ and content $\nu$.  This is straightforward: if the $i$th letter of $\ww$ is an $\hh$ (or $\vv$, respectively), then we extend the $i$th part of the $\ww$-decomposition of $\mu$ horizontally (or vertically, respectively) by $\nu_i^\ww$ boxes.  These new boxes are filled with the row numbers of the boxes in the $i$th part of the $\ww$-decomposition of $\nu$.  It is easy to verify that the above algorithm maximizes each $\pi_i^\ww$ in turn, so $\pi$ is $\ww$-maximal.
\qed\medskip

See Figure~\ref{fig_w} for an illustrated example of the proof.  We take $\ww=\hh\vv\vv\hh$.  In (a) and (b) we see the $\ww$-decompositions of $\mu=\pa{5411}$ with $\mu^\ww=(5,3,1,2)$ and of $\nu=\pa{3322}$ with $\nu^\ww=(3,3,3,1)$; each box in $\nu$ is labelled with its row.  In (c) we see the $\pi$ such that $\pi^\ww = (8,6,4,3) = \mu^\ww+\nu^\ww$; as a partition $\pi=\pa{8522211}$.  The boxes in $\pi/\mu$ are filled with the same numbers as the corresponding pieces of the $\ww$-decomposition of $\nu$ from (b), giving a Littlewood--Richardson tableau.

\begin{figure}
\begin{center}
\addtolength{\tabcolsep}{10pt}
\begin{tabular}{ccc}
%
\setlength{\unitlength}{17.5pt}
\begin{picture}(5,5)
\thicklines \hl045 \hl001 \hl123 \hl431 \vl004 \vl102 \vl421 \vl531 
\thinlines \hl034 \vl121 \vl221 
\end{picture}
&
\setlength{\unitlength}{17.5pt}
\begin{picture}(3,4)
\thicklines \hl043 \hl002 \hl221 \vl004 \vl202 \vl322 
\thinlines \hl033 \vl103 \vl221 
\multiput(0.5,3.5)(1,0){3}{\makebox(0,0){1}}
\multiput(0.5,2.5)(1,0){3}{\makebox(0,0){2}}
\multiput(0.5,1.5)(1,0){2}{\makebox(0,0){3}}
\multiput(0.5,0.5)(1,0){2}{\makebox(0,0){4}}
\end{picture}
&
\setlength{\unitlength}{10pt}  
\begin{picture}(8,7)(0,-3)
\thinlines 
\hl045 \hl035 \hl123 \hl001
\vl004 \vl103 \vl221 \vl421 \vl531
\thicklines 
\hl543 \hl434 \hl421 \hl121 \hl001 \hl1{-1}1 \hl0{-3}1
\vl0{-3}3 \vl1{-3}5 \vl2{-1}3 \vl421 \vl522 \vl831
{\scriptsize
\multiput(5,3)(1,0){3}{\makebox(1,1){1}}
\put(4,2){\makebox(1,1){2}}
\multiput(0,-1)(1,2){2}{\makebox(1,1){2}}
\multiput(0,-2)(1,2){2}{\makebox(1,1){3}}
\multiput(0,-3)(1,2){2}{\makebox(1,1){4}}}
\end{picture}
\\[6pt]
(a) $\mu=\pa{5411}$ & (b) $\nu=\pa{3322}$ & (c) $\pi=\pa{8522211}$
\end{tabular} \vspace{-1ex}
\end{center}
\caption{$\ww$-decompositions for $\ww=\hh\vv\vv\hh$, illustrating the lemma.}
\label{fig_w}
\end{figure}

If the word $\ww$ is $\hh\hh\ldots\hh$, the $\ww$-decomposition just breaks $\lambda$ into its rows, and the lemma reduces to the familiar fact that $s_{\mu+\nu}$ is the maximal term appearing in $s_\mu s_\nu$ in the dominance order.  A word of all $\vv$'s does the same for reverse dominance order, and other words somehow interpolate between the two.

\proof[Proof of Theorem]
Assume first that both dimensions of the rectangle $R$ are even; this is 
the case in which the proof is clearest.  Let $c=a/2$ and $d=b/2$ and pick 
any word $\ww$ with $c$ $\hh$'s and $d$ $\vv$'s, as above.  We will determine
the $\ww$-maximal partition $\pi$ that can appear in any product 
$s_\mu s_\muc$ for $\mu\subseteq R$.  By the lemma, we must maximize
$\mu_1^\ww+(\muc)_1^\ww$.  We do so by maximizing each of $\mu_1^\ww$
and $(\muc)_1^\ww$ independently: the $\ww$-decompositions of $\mu$ and $\muc$ 
each begin with an entire row or column of the rectangle (where $\ww$ begins with an
$\hh$ or a $\vv$, respectively).  Continuing inductively, we maximize
each $\mu_i^\ww$ and $(\muc)_i^\ww$ in turn; they are
independent and equal, because the $i$th pieces of the
corresponding $\ww$-decompositions are exchanged by rotating the
rectangle $180^\circ$.  We conclude that $\mu=\muc$.

Thus the $\ww$-maximal $s_\pi$ appears in $s_\lambda^2$,
where $\lambda$ is the $\ww$-maximal self-complementary partition,
and in no other $s_\mu s_\muc$.
The map from words $\ww$ to self-complementary $\lambda$ is a
bijection.  Each such $\lambda$ is determined by the path of length
$c+d$ running from the upper-right corner to the midpoint of the
rectangle, separating $\lambda$ from its complement.  The $i$th letter
in the corresponding word $\ww$ is $\hh$ or $\vv$ depending on whether
the $i$th step in this path is (confusingly) vertical or horizontal,
respectively.  See Figure~\ref{fig_46} for an illustration.

\begin{figure}
\begin{center}
\begin{tabular}{cccccc}
\begin{picture}(6,4)
\thicklines \hl006 \hl046 \vl004 \vl604 
\put(3,2){\circle*{.2}}
\thinlines \hl016 \hl026 \hl036 \vl104 \vl204 \vl304 \vl404 \vl504 
\end{picture}
&
\begin{picture}(6,4)
\thicklines \hl006 \hl046 \vl004 \vl604 
\put(3,2){\circle*{.2}}
\put(0,3){\makebox(1,1){$\scriptstyle\mu$}}
\put(5,0){\makebox(1,1){$\scriptstyle\muc$}}
\thicklines \hl016 \hl036 
\thinlines \hl026 \vl112 \vl212 \vl312 \vl412 \vl512 
\end{picture}
&
\begin{picture}(6,4)
\thicklines \hl006 \hl046 \vl004 \vl604 
\put(3,2){\circle*{.2}}
\put(0,3){\makebox(1,1){$\scriptstyle\mu$}}
\put(5,0){\makebox(1,1){$\scriptstyle\muc$}}
\thicklines \hl015 \hl135 \vl112 \vl512 
\thinlines \hl124 \vl112 \vl212 \vl312 \vl412 \vl512 
\end{picture}
&
\begin{picture}(6,4)
\thicklines \hl006 \hl046 \vl004 \vl604 
\put(3,2){\circle*{.2}}
\put(0,3){\makebox(1,1){$\scriptstyle\mu$}}
\put(5,0){\makebox(1,1){$\scriptstyle\muc$}}
\thicklines \hl014 \hl234 \vl212 \vl412 
\thinlines \hl222 \vl212 \vl312 \vl412 
\end{picture}
&
\begin{picture}(6,4)
\thinlines \hl006 \hl046 \vl004 \vl604 
\put(3,2){\circle*{.2}}\put(6,4){\circle*{.2}}
\put(0,3){\makebox(1,1){$\scriptstyle\mu$}}
\put(5,0){\makebox(1,1){$\scriptstyle\muc$}}
\thinlines \hl012 \hl221 \vl211 
\thicklines 
\put(6,4){\vector(0,-1){1}}
\put(6,3){\vector(-1,0){1}}
\put(5,3){\vector(-1,0){1}}
\put(4,3){\vector(0,-1){1}}
\put(4,2){\vector(-1,0){1}}
\end{picture}
\end{tabular}
\end{center}
\caption{Here $R$ is $4\times6$ and we pick $\ww=\hh\vv\vv\hh\vv$.  At each
stage we assign boxes of $R$ to $\mu$ or $\muc$ to make successive rows
or columns as long as possible, depending on whether the corresponding
letter in $\ww$ is an $\hh$ or $\vv$ respectively.  Note how the letters in $\ww$ 
give the $\mu$-$\muc$ boundary: 
$\hh =  
\protect\begin{picture}(.5,1)\thicklines
\protect\put(.25,1){\protect\vector(0,-1){1.25}} 
\protect\end{picture}$,
$\vv=
\protect\begin{picture}(1.5,1)\thicklines
\protect\put(1.25,.3){\protect\vector(-1,0){1.25}} 
\protect\end{picture}$.
}
\label{fig_46}
\end{figure}

Now assume that one dimension of $R$ is odd.  The above map from words $\ww$ 
(with $c=\frac{a-1}2$ or $d=\frac{b-1}2$) to self-complementary $\lambda$
still provides the proof. By the time the last step of $\ww$ that 
heads in the even-length direction is read, the maximizations fix the boundary 
between $\mu$ and $\muc$ and again they are equal.  Finally, if both dimensions 
of $R$ are odd and $\lambda$ is almost self-complementary, then the 
final maximization step leaves us free to assign the central square of $R$
to either $\mu$ or $\muc$.
\qed

\section{Conjecture}

A natural result that would subsume the theorem is the following:

\begin{conj}[weak version]
Fix a rectangle $R$.  The products $s_\lambda s_\lambdac$ are all
linearly independent.
\end{conj}

The set of products here is of course parametrized by {\em unordered}
pairs $(\lambda,\lambdac)$ complementary in $R$.  But the witnesses
that appear in the above proof show something much stronger than
linear independence, and it too has a natural extension:

\begin{conj}[strong version]
Fix a rectangle $R$, and define the sets:
\begin{eqnarray*}
P_0 &=& \{s_\lambda s_\lambdac\} \mbox{ as above}, \\
W_i &=& \{s_\lambda s_\lambdac \in P_i \,\,|\,\, \exists\pi : s_\lambda s_\lambdac \mbox{ is the only product in $P_i$  containing $s_\pi$}\},\\
P_{i+1} &=& P_i \setminus W_i.
\end{eqnarray*}
Then some $P_i$ is empty.
\end{conj}

In other words, if we take the set of all products $s_\lambda s_\lambdac$ 
and iterate the operation ``eliminate any products which can now be proved 
linearly independent by some witness $s_\pi$,'' then we eventually eliminate
all products.  This is equivalent to the statement that with appropriate orderings 
of the set of products and Schur functions, the matrix expanding the products in 
the Schur basis is upper triangular.

Note an interesting difference between the weak and
strong versions of the conjecture.  Verifying linear independence as
in the weak version requires, a priori, knowing the actual values of
all of the Littlewood--Richardson coefficients
$c_{\lambda\lambdac}^\nu$.  The strong version only makes claims about
which of these coefficients are {\em nonzero}.  The collection of
triples of partitions $\{(\lambda,\mu,\nu) \,|\, c_{\lambda\mu}^\nu
\neq 0\}$ consists of all integer points in some cone, according to
the Saturation Conjecture, now a theorem of Knutson and Tao~\cite{honey1}; 
moreover the minimal set of linear inequalities defining the cone is now 
known~\cite{honey2}.  The strong version of the conjecture, if true, follows 
from this set of linear inequalities.

The remainder of this section gives evidence supporting this
conjecture: computer verification for small rectangles, proofs in some
special cases, and the beginning of a sadly incomplete inductive proof.
Finally we mention a few failed attempts at generalization.

\subsection{Computer observations}

The strong version of the conjecture has been verified by computer, using
Anders Buch's Littlewood--Richardson calculator,  for all $R$ fitting inside 
the $8\by8$ or $7\by9$ rectangle.  For a sense of scale, the $8\by8$
calculation involves 6470 products, and they are linear combinations of 
395,377 distinct Schur functions (out of the 487,842 Schur functions of 
degree 64 whose partitions fit in a $16\by16$ box).


The computer results confirm that the conjecture still holds if we require the
witness $s_\pi$ to always have coefficient one.  That is, with respect to
some ordering, the matrix writing the products in the Schur basis 
appears to be upper unitriangular. 

The theorem is the statement that the set $W_0$ contains all of the
self-complementary or almost self-complementary products.  According
to the computer calculations, these appear to be the entirety of $W_0$.

\subsection{Quick Proofs for $1\by b$ and $2\by b$}

When $R$ is a single row (or column), the conjecture is checked
easily, e.g.:
$$
\setlength{\unitlength}{8pt}
\begin{array}{lcl}
\rr30 \rr30 &=& \rrr600 + \rr51 + \rr42 + \rr33 \\
\rr40 \rr20 &=& \rrr600 + \rr51 + \rr42  \\
\rr50 \rr10 &=& \rrr600 + \rr51  \\
\rr60 \rr00 &=& \rr60 \\
\end{array}
$$
Here we used that $s_\pa{j} s_\pa{k} = \sum s_\pa{j+k-i,i}$ for
$i=0,1,\ldots,\min(j,k)$.

We can also give explicit witnesses for the $2\by b$ case.  To fix
notation, let us always choose $\lambda=\pa{\lambda_1,\lambda_2}$ to
be larger than $\lambdac$, so that $\lambda_1 + \lambda_2 \geq b$.  We
partially order the products $s_\lambda s_\lambdac$ by increasing
values of $\lambda_1+\lambda_2$; note that this ordering begins with
the cases $\lambda_1+\lambda_2=b$, which holds precisely when
$\lambda=\lambdac$.  The general witness for $\lambda$ is
$\pi=\pa{2\lambda_1-\lambda_2, \lambda_2, b-\lambda_1, b-\lambda_1}$.
There are no choices in constructing the unique Littlewood--Richardson
tableau of shape $\pi/\lambda$ and content $\lambdac$: 
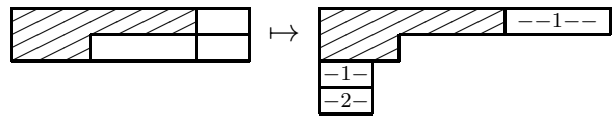
$$
\begin{picture}(9,2)
\thicklines \hl003 \hl314 \hl027 \vl002 \vl301 \vl711
\thinlines
\multiput(0,0)(1,0){3}{\line(2,1){1}}
\multiput(0,.5)(1,0){3}{\line(2,1){1}}
\multiput(0,1)(1,0){7}{\line(2,1){1}}
\multiput(0,1.5)(1,0){7}{\line(2,1){1}}
\thicklines \hl306 \hl712 \hl722 \vl701 \vl902
\end{picture}
\,\, \rb{.75}{${}\mapsto{}$} \,\,
\rb{-2}{\begin{picture}(11,4)(0,-2)
\thicklines \hl003 \hl314 \hl027 \vl002 \vl301 \vl711
\thinlines
\multiput(0,0)(1,0){3}{\line(2,1){1}}
\multiput(0,.5)(1,0){3}{\line(2,1){1}}
\multiput(0,1)(1,0){7}{\line(2,1){1}}
\multiput(0,1.5)(1,0){7}{\line(2,1){1}}
\thicklines \hl714 \hl724 \vl{11}11
\hl0{-1}2 \hl0{-2}2 \vl0{-2}2 \vl2{-2}2
\put(0,-1){\makebox(2,1){$\scriptstyle{}-1-{}$}}
\put(0,-2){\makebox(2,1){$\scriptstyle{}-2-{}$}}
\put(7,1){\makebox(4,1){$\scriptstyle{}--1--{}$}}
\end{picture}}
$$
To see that $s_\pi$ has coefficient zero in all products larger in the
ordering, observe that the fourth row of the witness forces
$\lambda_1$ to be small, but its first row forces $\lambda_1 +
(\lambdac)_1$ to be as large as possible.  This argument is left
sketchy in part because we will see one that subsumes it shortly.

\subsection{A Little Induction}

We can find some witnesses by induction from smaller rectangles.

Suppose $R$ is $a\by b$ and we have a complementary pair $\lambda$,
$\lambdac$ such that the entire top row of $R$ is in $\lambda$ and the
entire bottom row is in $\lambdac$, i.e. $\lambda_1 = (\lambdac)_1 =
b$.  Then strip off these two rows to get a pair $\overline{\lambda}$,
$\overline{\lambdac}$ in the $(a-2)\by b$ rectangle. 

Suppose their product has a witness $\overline{\pi}$.  Then $\pi$ is a
witness for $\lambda$, $\lambdac$, where $\pi$ is constructed by
inserting a row of size $2b$ before the first row of $\overline{\pi}$:
The only complementary products $s_\mu s_\muc$ that could possibly
contain any $s_\pi$ with $\pi_1=2b$ are those in which
$\mu_1=(\muc)_1=b$, and $\pi$ appears in such products if and only if
the associated $\overline{\pi}$ appears in the corresponding product
$s_{\overline{\mu}} s_{\overline{\muc}}$. 

Likewise, if $\lambda$ contains the leftmost column of $R$ and
$\lambdac$ contains the rightmost, we can bootstrap a witness obtained
from the $a\by (b-2)$ rectangle by adding a column of height $2a$.  It
seems natural to hope for another induction for the remaining case, in
which $\lambda$ contains both the top row and left column, but I cannot
provide one.

This limited induction argument yields a witness-based proof of linear
independence for the $(\lambda, \lambdac)$ as long as the recursion
never gets to a case where $\lambda$ contains both the top row and left
column.  But such $(\lambda, \lambdac)$ are precisely the
self-complementary pairs of the previous section; the witness
bootstrapping process here gives the same witnesses constructed in the
proof of the theorem.

This also gives a simplified proof of the $2\by b$ case: by the
induction step, we need only consider the situation when
$\lambda=\pa{b,\lambda_2}$, in which case the sequence of witnesses
$\pi=\pa{2b-\lambda_2, \lambda_2}$ is transparent.  Again, this
bootstraps to the same witness already given explicitly.

\subsection{Generalizations?}

To save the interested reader some time, we mention a few
plausible-sounding generalizations of the conjecture which appear
to fail.

Complementation with respect to a rectangle is a natural operation
from the Grassmannian point of view.  But there and in other similar
settings one takes restricted products of Schur functions, by
specializing $s_\nu\mapsto0$ if $\nu$ has too many rows or columns or
both.  These tend to fail almost immediately (e.g. the $2\by2$ case in
the Introduction) just by running out of room.

Given an arbitrary partition $R$, there are two natural sets of
products which specialize to $\{s_\lambda s_\lambdac\}$ when $R$ is a
rectangle.  First, we could consider the set $\{s_\lambda
s_{R/\lambda}\}$ for all $\lambda\subseteq R$; the special case of
rectangular $R$ is the only one where the skew Schur function
$s_{R/\lambda}$ is again a single Schur function.  But for $R=\pa{21}$
we have the linear dependence $s_\pa{1} s_{R/\pa{1}} = s_\pa{2} s_{R/\pa{2}} + 
s_\pa{11} s_{R/\pa{11}}$.

Second, we could consider the set of products $\{s_\lambda s_\mu\}$ in
which a term $s_R$ appears; this is closely related to the coproduct
of the Hopf algebra structure on the ring of symmetric functions.  A
counterexample here comes from considering the two ways to
parenthesize the product $\Inline{\rr20\rr10\rr10}$.  We find that
$(\Inline{\rr30}+\Inline{\rr21})\Inline{\rr10} =
\Inline{\rr20}(\Inline{\rr20}+\Inline{\rr11})$, yet all four products
contain a term $\Inline{\rr31}$.

Finally, we could hope to replace Schur functions with their type-B analogues,
the Schur-Q functions (or Schur-P functions, which differ only by a scalar, 
irrelevant here).  They are indexed by strict partitions, often represented 
graphically by shifted diagrams.  The role of the rectangle is played by the
staircase partitions $S=\pa{n\ldots321}$, and the complement of 
$\lambda\subseteq S$ is the remaining boxes of $S$ reflected through
the line of slope 1.  For example, $\inline{\Srr41}$ and $\inline{\Srr32}$ are
complementary in $\inline{{\rb{-3}{{\yd{(4,4)(-4,0)}{\Srow4{4}\Srow3{3}\Srow2{2}\Srow1{1}
\thicklines \vl{-2}21 \hl{-2}32}}}}}$.  The conjecture fails for $S=\pa{321}$, as
$
\Inline{\Srrr321} + \Inline{\Srr31\Srr20} = \Inline{\Srr32\Srr10} + \Inline{\Srr30\Srr21}.
$



\begin{thebibliography}{9}

\bibitem{adler}
M. Adler, V. B. Kuznetsov, and P. van Moerbeke.
Rational solutions to the Pfaff lattice and Jack polynomials.
Preprint nlin.SI/0202037 (February 2002).

\bibitem{honey1}
A. Knutson, T. Tao.
The honeycomb model of ${\rm GL}_n$ tensor products I: Proof of the
saturation conjecture. 
{\em J. Amer. Math. Soc.}  {\bf 12}  (1999),  no. 4, 1055--1090.

\bibitem{honey2}
A. Knutson, T. Tao, C. Woodward.
The honeycomb model of ${\rm GL}_n$ tensor products II: Puzzles
determine facets of the Littlewood-Richardson cone. 
{\em J. Amer. Math. Soc.}, to appear.

\end{thebibliography}
\end{document}